\documentclass[12pt]{amsart}
\usepackage{amssymb,latexsym,amsthm}
\usepackage{graphicx}
\usepackage{color}
\usepackage{mathrsfs}
\usepackage[top=30truemm, bottom=30truemm, left=25truemm, right=25truemm]{geometry}
\usepackage{amsmath,amssymb}

\newtheorem{theorem}{Theorem}
\newtheorem{lemma}[theorem]{Lemma}
\newtheorem{cor}[theorem]{Corollary}
\newtheorem{prop}[theorem]{Proposition}

\theoremstyle{definition}
\newtheorem{definition}[theorem]{Definition}
\newtheorem{example}[theorem]{Example}

\theoremstyle{remark}

\newcommand{\1}{\mathbf{1}}

\newcommand{\ch}{\operatorname{Ch}}
\newcommand{\M}{\mathfrak M}

\newcommand{\nn}{\|\hspace{-0.45mm}|}

\numberwithin{equation}{section}

\newcommand{\bit}{\begin{itemize}}
\newcommand{\eit}{\end{itemize}}

\newcommand{\bea}{\begin{eqnarray*}}
\newcommand{\eea}{\end{eqnarray*}}

\newcommand{\wbo}{\widetilde{B_1}}
\newcommand{\wbt}{\widetilde{B_2}}
\newcommand{\wbj}{\widetilde{B_j}}
\newcommand{\wb}{\widetilde{B}}

\newcommand{\Lip}{\operatorname{Lip}}
\newcommand{\lip}{\operatorname{lip}_{\alpha}}

\makeatletter 
\newcommand{\emsideset}[3]{%
  \@mathmeasure\z@\displaystyle{#3}%
  \global\setbox\@ne\vbox to\ht\z@{}\dp\@ne\dp\z@ 
  \setbox\tw@\box\@ne 
  \@mathmeasure4\displaystyle{\copy\tw@#1}%
  \@mathmeasure6\displaystyle{#3#2}%
  \dimen@-\wd6 \advance\dimen@\wd4 \advance\dimen@\wd\z@ 
  \hbox to\dimen@{}\mathop{\kern-\dimen@\box4\box6}%
} 
\makeatother 

\begin{document}
\title[Isometries]
{Isometries on Banach algebras of $C(Y)$-valued maps
}

\author{
Osamu~Hatori
}
\address{
Department of Mathematics, Faculty of Science,
Niigata University, Niigata 950-2181, Japan
}
\email{hatori@math.sc.niigata-u.ac.jp
}






\begin{abstract}
We propose a unified approach to the study of isometries on algebras of vector-valued Lipschitz maps and those of continuously differentiable maps by means of the notion of natural $C(Y)$-valuezations that take values in unital commutative $C^*$-algebras. 
\end{abstract}

\maketitle

\section{Introduction}
The study on isometries on Banach algebras dates back to the classical Banach-Stone theorem. After that there are many literature on the study of isometries not only on Banach algebras but also Banach spaces of functions and operators. In this paper we study isometries on the algebra of Lipschitz functions and continuously differentiable functions with values in unital commutative $C^*$-algebras. We propose a unified approach to such a study by considering natural $C(Y)$-valuezations.

The study on the space of Lipschitz functions is probably initiated by
de Leeuw \cite{dl} for functions on the real line. Roy \cite{roy} considered isometries on the Banach space $\Lip(K)$ of Lipschitz functions on a compact metric space $K$, equipped with the norm $\|f\|_M=\max\{\|f\|_\infty, L(f)\}$, where $L(f)$ denotes the Lipschitz constant. On the other hand, 
Cambern \cite{c} studied isometries on spaces of continuously differentiable functions $C^1([0,1])$ with norm given by $\|f\|=\max_{x\in [0,1]}\{|f(x)|+|f'(x)|\}$ for $f\in C^1([0,1])$ and exhibited the forms of the surjective isometries supported by such spaces. 
Rao and Roy \cite{rr} proved that surjective isometries between $\Lip([0,1])$ and $C^1([0,1])$ with respect to the norm $\|f\|_L=\|f\|_\infty+\|f'\|_{\infty}$ are of the canonical forms in the sense that they are weighted composition operators. Jim\'enez-Vargas and Villegas-Vallecillos in \cite{amPAMS} considered isometries of spaces of vector-valued Lipschitz maps on a compact metric space taking values in a strictly convex Banach space, equipped with the norm $\|f\|=\max\{\|f\|_\infty, L(f)\}$, see also \cite{amHouston}. Botelho and Jamison \cite{bjStudia2009} studied isometries on $C^1([0,1],E)$ with $\max_{x\in [0,1]}\{\|f(x)\|_E+\|f'(x)\|_E\}$.
See also \cite{mw,amy,araduba,bfj,kos,bjz,rm,bjPositivity17,mt,kawar,kc1,kc12,lcmw,kkm,lpww,jlp}. Refer also a book of Weaver \cite{weab}.

In this paper an isometry means a complex-linear isometry.
 Isometries on algebras of Lipschitz maps and continuously differentiable maps have often been studied independently. 
Jarosz \cite{ja} and Jarosz and Pathak \cite{jp} studied a problem when an isometry on a space of continuous functions is a weighted composition operator. They provided a unified approach for function spaces such as $C^1(K)$, $\Lip(K)$, $\lip(K)$ and $AC[0,1]$. In particular, Jarosz \cite[Theorem]{ja} proved that a \emph{unital} surjective isometry between unital semisimple commutative Banach algebras with \emph{natural norms} is canonical. 

We consider a Banach algebra of continuous maps defined on a compact Hausdorff space whose values are in a unital $C^*$-algebra. It is an abstraction of $\Lip(K,C(Y))$ and $C^1(K,C(Y))$.  
We propose a unified approach to the study of isometries on algebras  $\Lip(K,C(Y))$, $\lip(K,C(Y))$ and $C^1(K,C(Y))$, where $K$ is a compact metric space, $[0,1]$ or $\mathbb{T}$ (in this paper $\mathbb{T}$ denotes the unit circle on the complex plane),  and $Y$ is a compact Hausdorff space. We study isometries without assuming that they preserve unit. We prove that the form of isometries between such algebras are of a canonical form. 
As corollaries of the result, we describe isometries on $\Lip(K,C(Y))$, $\lip(K,C(Y))$, $C^1([0,1],C(Y))$, and $C^1(\mathbb{T},C(Y))$ respectively.

The main result Theorem 14 in \cite{hatori} with a detailed proof is recaptured as Theorem \ref{main} in this paper. It gives the form of a surjective isometry $U$ between certain Banach algebra of continuous maps with values in unital $C^*$-algebras. 
Verifying  that $U(1)=1\otimes h$ for an $h\in C(Y_2)$ with $|h|=1$ on $Y_2$ due to Choquet's theory (Proposition \ref{absolute value 1}), 
the Lumer's method (cf. \cite{hoacta}) works very well. We see that $U_0$ is a composition operator of type BJ (cf. \cite{hot}). Then \ref{main} is proved.

This paper surveys recent papers \cite{hoacta} by Hatori and Oi, and \cite{hatori} by Hatori.

\section{Preliminaries}
Let $Y$ be a compact Hausdorff space and $E$ a real or complex Banach space. The space of all $E$-valued  continuous maps on $Y$ is denoted by $C(Y,E)$. When $E={\mathbb C}$ (resp. $\mathbb{R}$), $C(Y,E)$  is  abbreviated by $C(Y)$ (resp. $C_{\mathbb R}(Y)$). 
The supremum norm on $S\subset Y$ is  $\|F\|_{\infty(S)}=\sup_{x\in S}\|F(x)\|_E$ for $F\in C(Y,E)$. 
We may omit the subscript $S$ and write only $\|\cdot\|_{\infty}$. 
Let $K$ be a compact metric space and $0<\alpha\le 1$. Put 
\[
L_\alpha(F)=\sup_{x\ne y}\frac{\|F(x)-F(y)\|_E}{d(x,y)^\alpha}.
\]
for $F\in C(K,E)$. The number $L_\alpha$ is called the $\alpha$-Lipschitz number of $F$. When $\alpha=1$ we  omit the subscript $\alpha$ and write only $L(F)$ and call it the Lipschitz number. We denote 
\[
\Lip_\alpha(K,E)=\{F\in C(K,E):L_\alpha(F)<\infty\}.
\]
When $E={\mathbb C}$, $\Lip(K,{\mathbb C})$ is abbreviated to $\Lip(K)$.
When $\alpha=1$ the subscript is omitted and it is written as $\Lip(K,E)$ and $\Lip(K)$. 
When $0<\alpha<1$ the subspace
\begin{equation*}
\lip(K,E)
=\{F\in \Lip_\alpha(K,E):\text{$\lim_{x\to x_0}\frac{\|f(x_0)-f(x)\|_E}{d(x_0,x)^\alpha}=0$ for every $x_0\in K$}\}
\end{equation*}
of $\Lip_\alpha(K,E)$ is called a little Lipschitz space. For $E={\mathbb C}$, $\lip(K,{\mathbb C})$) is abbreviated to $\lip(K)$. 
A variety of complete norms on $\Lip_\alpha(K,E)$ and $\lip(K,E)$ exist. 
The norm $\|\cdot\|_L$ of $\Lip_\alpha(K,E)$ (resp. $\lip(K,E)$) is defined by 
\[
\|F\|_L=\|F\|_{\infty(K)}+L_\alpha(F), \quad F\in \Lip_\alpha(K,E)\,\, \text{(resp. $\lip(K,E)$)},
\]
which is often called the $\ell^1$-norm or the sum norm.
The norm, which is called the max norm, $\|\cdot\|_M$ of $\Lip_\alpha(K,E)$ (resp. $\lip(K,E)$) is defined by 
\[
\|F\|_M=\max\{\|F\|_\infty,\,\, L_\alpha(F)\},\quad F\in \Lip_\alpha(K,E)\,\, \text{(resp. $\lip(K,E)$)}.
\]
Note that $\Lip_{\alpha}(K,E)$ (resp. $\lip(K,E)$) is a Banach space with respect to $\|\cdot\|_L$ and $\|\cdot\|_M$ respectively. If $E$ is a unital (commutative) Banach algebra, then the norm $\|\cdot\|_L$ is sub-multiplicative. Hence $\Lip_{\alpha}(K,E)$ (resp. $\lip(K,E)$) is a unital (commutative) Banach algebra with respect to the norm $\|\cdot\|_L$ if $E$ a unital (commutative) Banach algebra. The norm $\|\cdot\|_M$ needs not be sub-multiplicative even if $E$ is a Banach algebra. With the norm $\|\cdot\|_M$, $\Lip_{\alpha}(K)$ and $\lip(K)$ need not be Banach algebras . 
In this paper we mainly concerns with $\|\cdot\|_L$ 
and $E=C(Y)$. Then $\Lip_\alpha(K,C(Y))$ and $\lip(K,C(Y))$ are unital semisimple commutative Banach algebras with $\|\cdot\|_L$.  

Let $K=[0,1]$ or ${\mathbb T}$. 
We say that $F\in C(K,E)$ is continuously differentiable if there exists $G\in C(K,E)$ such that
\[
\lim_{K\ni t\to t_0}\left\|\frac{F(t_0)-F(t)}{t_0-t}-G(t_0)\right\|_E=0
\]
for every $t_0\in K$. We denote $F'=G$. 
Put 
\[
C^1(K,E)=\{F\in C(K,E):\text{$F$ is continuously differentiable}\}.
\]
Then $C^1(K,E)$ with norm $\|F\|=\|F\|_\infty+\|F'\|_\infty$ is a Banach space and it is unital (commutative) Banach algebra provided that $E$ is a unital (commutative) Banach algebra. We mainly consider the case where $E=C(Y)$ with the supremum norm for a compact Hausdorff space $Y$. In this case  $C^1(K,C(Y))$ with the norm $\|F\|=\|F\|_\infty+\|F'\|_\infty$ for $F\in C^1(K,C(Y))$ is a unital semisimple commutative Banach algebra. 
We may suppose that $C(Y)$ is isometrically isomorphic to ${\mathbb C}$ if $Y$ is a singleton, and we abbreviate $C^1(K,C(Y))$ by $C^1(K)$ when $Y$ is a singleton.

By identifying $C(K, C(Y))$ with $C(K\times Y)$ we may assume that $\Lip(K,C(Y))$ (resp. $\lip(K,C(Y))$) is a subalgebra of $C(K\times Y)$ by the correspondence 
\[
F\in \Lip(K,C(Y)) \leftrightarrow ((x,y)\mapsto (F(x))(y))\in C(K\times Y).
\]
Under this identification we may suppose that 
$\Lip(X,C(Y))\subset C(X\times Y)$, 
$\lip(X,C(Y))\subset C(X\times Y)$, and
$C^1(K,C(Y))\subset C(K\times Y)$.

Let $\emptyset \ne Q\subset C(Y)$. We say that $Q$ is point separating or $Q$ separates the points of $Y$ 
if 
for every pair $x$ and $y$ of distinct points in $Y$, there corresponds a function $f\in Q$ such that $f(x)\ne f(y)$.
 In this paper, unity of a Banach algebra $B$ is denoted by $\1$. The maximal ideal space of $B$ is denoted by $M_B$.
\section{A theorem of Jarosz on isometries which preserve $\1$}
In most cases the form of an isometry between Banach algebras depends not only on the algebraic structure, but also on the norms on theses algebras.
Jarosz \cite{ja} introduced \emph{natural norms} on spaces of continuous functions. He proved that isometries between a variety of spaces of continuous functions equipped with the natural norms are of canonical forms. 
See \cite{ja} for precise notations and terminologies. The following is a theorem of Jarosz on surjective unital isometries \cite{ja}.
\begin{theorem}[Jarosz \cite{ja}]\label{ja}
Let $X$ and $Y$ be compact Hausdorff spaces, let $A$ and $B$ be complex linear subspaces of $C(X)$ and $C(Y)$, respectively, and let $p,q\in \mathcal{P}$. Assume $A$ and $B$ contain constant functions, and let $\|\cdot\|_A$, $\|\cdot\|_B$ be a $p$-norm and $q$-norm on $A$ and $B$, respectively. Assume next that there is a linear isometry $T$ from $(A,\|\cdot\|_A)$ onto $(B,\|\cdot\|_B)$ with $T\1=\1$. Then if $D(p)=D(q)=0$, or if $A$ and $B$ are regular subspaces of $C(X)$ and $C(Y)$, respectively, then $T$ is an isometry from $(A,\|\cdot\|_\infty)$ onto $(B,\|\cdot\|_\infty)$.
\end{theorem}
We provide a precise proof of a theorem of Jarosz in \cite{hatori} by making an ambitious revision of one in \cite{ja}. 

In the following  a unital semisimple commutative Banach algebra $A$ is identified through the Gelfand transform with a unital subalgebra of $C(M_A)$ for maximal ideal space $M_A$ of $A$. Hence we see that the uniform closure of a unital semisimple commutative Banach algebra in $C(M_A)$ is a uniform algebra on $M_A$.
A unital semisimple commutative Banach algebra is regular in the sense of Jarosz \cite{ja}. Applying a theorem of Nagasawa \cite{naga} (cf. \cite{drw}) we have the following.
\begin{cor}[Corollary 2 \cite{hatori}]\label{semisimplecase}
Let $A$ and $B$ be unital semisimple commutative Banach algebras with natural norms. 
Suppose that $T:A\to B$ is a surjective complex-linear isometry with $T\1=\1$. Then there exists a homeomorphism $\varphi:M_B\to M_A$ such that 
\[
T(f)(x)=f\circ\varphi(x),\quad f\in A,\,\,x\in M_B.
\]
In particular, $T$ is an algebra isomorphism.
\end{cor}
We omit a proof (see the proof of Corollary 2 in \cite{hatori}).
\begin{cor}[Corollary 3 \cite{hatori}]\label{1lip}
Let $K_j$ be a compact metric space for $j=1,2$. Suppose that $T:\Lip(K_1)\to \Lip(K_2)$ is a surjective complex-linear isometry with respect to the norm $\|\cdot\|_L$. Assume $T\1=\1$. Then there exists a surjective isometry $\varphi:X_2\to X_1$ such that 
\begin{equation}\label{lipkl}
Tf(x)=f\circ \varphi(x),\quad f\in \Lip(K_1),\,\,x\in K_2.
\end{equation}
Conversely if $T:\Lip(K_1)\to \Lip(K_2)$ is of the form as \eqref{lipkl}, then $T$ is a surjective isometry with respect to both of $\|\cdot\|_M$ and $\|\cdot\|_L$ such that $T\1=\1$.
\end{cor}
We exhibit a proof which is a little bit precise than one given in \cite{hatori}.
\begin{proof}
It is well known that $(\Lip(K_j),\|\cdot\|_L)$ is a unital semisimple commutative Banach algebra with maximal ideal space $K_j$. Hence  Corollary \ref{semisimplecase} asserts that there is a homeomorphism $\varphi:K_2\to K_1$ such that 
\begin{equation}\label{HH1}
Tf(x)=f\circ \varphi(x),\quad f\in \Lip(K_1),\,\,x\in K_2.
\end{equation}
Let $y_0\in K_2$. Define $f:K_1\to {\mathbb C}$ by $f(x)=d_1(x,\varphi(y_0))$. Then by a simple calculation we infer that $L(f)=1$. Since $T$ is an isometry with respect $\|\cdot\|_L$, so is for $\|\cdot\|_\infty$ by Corollary \ref{semisimplecase}. Hence $L(Tf)=1$. By the definition of $L(\cdot)$ we have
\[
1=L(Tf)\ge \frac{d_1(\varphi(x),\varphi(y_0)}{d_2(x,y_0)}
\]
for every $x\in K_2\setminus \{y_0\}$. Thus we have that $d_2(x,y_0)\ge d_1(\varphi(x),\varphi(y_0))$ for every $x\in K_2$. As $y_0\in K_2$ is arbitrary, we have that 
\begin{equation}\label{d2d1}
d_2(x,y)\ge d_1(\varphi(x),\varphi(y))
\end{equation}
for every pair $x$ and $y$ in $K_2$. By \eqref{HH1} we have 
\[
T^{-1}g(x)=g\circ \varphi^{-1}(x), \quad g\in \Lip(K_2), \,\, x\in K_1.
\]
Since $T^{-1}$ is an isometry we have the same argument as above that
\begin{equation}\label{d1d2}
d_1(x,y)\ge d_2(\varphi^{-1}(x),\varphi^{-1}(y))
\end{equation}
for every pair $x$ and $y$ in $K_1$. By \eqref{d2d1} and \eqref{d1d2} we have that $\varphi$ is an isometry.

We omit a proof of the converse statement since it is trivial.
\end{proof}
It is natural to ask what is the form of a surjective isometry between $\Lip(K)$ without the hypothesis of $T\1=\1$.
Rao and Roy \cite{rr} proved that it is a weighted composition operator if $K=[0,1]$. They asked whether a surjective isometry on $\Lip(K)$ with respect to the metric induced by $\ell^1$-norm is induced by an isometry on $K$. Jarosz and Pathak \cite[Example 8]{jp} exhibited a positive solution. After the publication of \cite{jp} some authors expressed their suspicion about the argument there and the validity of the statement there had not been confirmed. In \cite{hoacta} we proved that Example 8 in \cite{jp} is true. In this paper we exhibit a slight general result (see also \cite{hatori}).

Since the max norm $\|\cdot\|_M$ is not sub-mutiplicative in general, 
$(\Lip(K),\|\cdot\|_M)$ need  not be a Banach algebra. By a simple calculation it is easy to see that $\|\cdot\|_M$ is a natural norm in the sense of Jarosz (see \cite{ja}) such that $\lim_{t\to +0}\frac{\max\{1,t\}-1}{t}=0$. By Theorem \ref{ja} we have the following. Refer the proof of Corollary 4 in \cite{hatori}.
\begin{cor}[Corollary 4 \cite{hatori}]\label{2lip}
Let $K_j$ be a compact metric space for $j=1,2$. Suppose that $T:\Lip(K_1)\to \Lip(K_2)$ is a surjective complex-linear isometry with respect to the norm $\|\cdot\|_M$. Assume $T\1=\1$. Then there exists a surjective isometry $\varphi:X_2\to X_1$ such that 
\begin{equation}\label{lipkm}
Tf=f\circ \varphi,\qquad f\in \Lip(K_1).
\end{equation}
Conversely if $T:\Lip(K_1)\to \Lip(K_2)$ is of the form as \eqref{lipkm}, then $T$ is a surjective isometry with respect to both of $\|\cdot\|_M$ and $\|\cdot\|_L$ such that $T\1=\1$.
\end{cor}
When 
 $T\1=\1$ is not assumed in Corollary \ref{2lip}, a simple counterexample such that $K_j$ is a two-point-set is given by Weaver\cite[p.242]{wea} (see also \cite{weab}) shows that $T$ need not be a weighted composition operator.

We have already pointed out \cite{hjv} that the original proof of Theorem \ref{ja} need a revison and made an ambitious revision in \cite{hjv,hoacta}. Although the revised proof  for a general case \cite{hoacta} is similar to that of Proposition 7 in \cite{hjv}, a detailed revision is exhibited in \cite{hatori}. Note that Tanabe \cite{tanabe} pointed out that $\lim_{t\to +0}(p(1,t)-1)/t$ always exists and it is finite for every $p$-norm.
To prove Theorem \ref{ja} we need Lemma 2 in \cite{ja} in the same way as the original proof of Jarosz. We note minor points in the original proof of Lemma 2. 
Note first that  five $\varepsilon/2$'s between 11 lines and 5 lines from the bottom of page 69 read as $\varepsilon/3$. Next $x\in X\setminus U_1$ reads as $x\in U_1$ on the bottom of page 69.  
We point out that the term $\displaystyle{\sum_{j=1}^{k_0-1}(f_j(x)-1)}$ which appears on the first line of the first displayed inequalities on page 70 reads $0$ if $k_0=1$. 
The term $1+\varepsilon$ on the right hand side of the second line of the same inequalities reads as $1+\frac{\varepsilon}{3}$. Two $\frac{\varepsilon}{2}$'s on the same line read as $\frac{\varepsilon}{3}$. On the next line $\frac{n+1}{n}\frac{\varepsilon}{2}$ reads as $\frac{\varepsilon}{3}$. 
For any $1\le k_0\le n$ we infer that 
\[
1\ge 1-2\frac{k_0-1}{n}\ge 1-2\frac{n-1}{n}>-1.
\]
Hence we have $|f(x)|\le 1+\varepsilon$ if $x\in U_1$ by the first displayed inequalities of page 70. 
The inequality $\|f\|_{\infty}\le \varepsilon$ 
on the fifth line on page 70 reads as $\|f\|_{\infty}\le1+\varepsilon$.

A precise revision of the proof of a theorem of Jarosz (Theorem \ref{ja}) is given in \cite{hatori}.
\section{Banach algebras of $C(Y)$-valued maps}
Let $X$ be a compact Hausdorff space and
 $B$  a unital point separating subalgebra of $C(X)$ equipped with a Banach algebra norm. Then $B$ is semisimple since  
$\{f\in B:f(x)=0\}$ is a maximal ideal of $B$ for every $x\in X$ and the Jacobson radical of $B$ vanishes. The inequality $\|f\|_\infty\le \|f\|_B$ for every $f\in B$ is well known. 
We say that $B$ is \emph{natural} Banach algebra if the map $e:Y\to M_B$ defined by $y\mapsto \phi_y$, where $\phi_y(f)=f(y)$ for every $f\in B$, is bijective. 
We say that $B$ is \emph{self-adjoint} if $B$ is natural and conjugate-closed in the sense that $f\in B$ implies that $\bar{f}\in B$ for every $f\in B$, where $\bar{\cdot}$ denotes the complex conjugation on $Y$. 

Let $X$ and $Y$ be compact Hausdorff spaces. 
For functions $f\in C(X)$ and $g\in C(Y)$, let $f\otimes g\in C(X\times Y)$ be the function defined by $f\otimes g(x,y)=f(x)g(y)$ for $(x,y)\in X\times Y$, and for a subspace $E_X$ of $C(X)$ and a subspace $E_Y$ of $C(Y)$, let
\[
E_X\otimes E_Y=\left\{\sum_{j=1}^nf_j\otimes g_j: n\in {\mathbb N},\,\,f_j\in E_X,\,\,g_j\in E_Y\right\},
\]
and
\[
\1\otimes E_Y=\{\1\otimes g: g\in E_Y\}.
\]
A natural $C(Y)$-valuezation is introduced in \cite{hatori}.
\begin{definition}[Definition 12 in \cite{hatori}]\label{n}
Let $X$ and $Y$ be compact Hausdorff spaces. Suppose that $B$ is a unital point separating subalgebra of $C(X)$ equipped with a Banach algebra norm $\|\cdot\|_B$. Suppose that $B$ is self-adjoint. Suppose that $\wb$ is a unital point separating subalgebra of $C(X\times Y)$ such that $B\otimes C(Y)\subset \wb$ equipped with a Banach algebra norm $\|\cdot\|_{\wb}$.  Suppose that $\wb$ is self-adjoint. 
We say that $\wb$ is a \emph{natural $C(Y)$-valuezation of $B$} if
there exists a compact Hausdorff space $\M$ and a complex-linear map $D:\wb\to C(\M)$ such that $\ker D=1\otimes C(Y)$ and $D(C_{\mathbb R}(X\times Y)\cap \wb)\subset C_{\mathbb R}(\M)$ which satisfies 
\[
\|F\|_{\wb}=\|F\|_{\infty(X\times Y)}+\|D(F)\|_{\infty(\M)},\quad F\in \wb.
\]
\end{definition}
 Note that the norm $\|\cdot\|_{\wb}$ is a natural norm in the sense of Jarosz \cite{ja}. 

Let $B$ be a Banach algebra which is a unital separating subalgebra of $C(X)$ for a compact Hausdorff space $X$. Then if $(X,C(Y),B,\wb)$ is an admissible quadruple of type L defined in \cite{hoacta}, then $\wb$ is a natural $C(Y)$-valuezation of $B$ due to Definition  \ref{n}. On the other hand for a $C(Y)$-valuezation of $\wb$ of $B$, $(X,C(Y),B,\wb)$ need not be an admissible quadruple defined by Nikou and O'Farrell \cite{no} (cf. \cite{hot}). This is because we do not assume that $\{F(\cdot,y):F\in \wb,\,\,y\in Y\}\subset B$, which is a requirement for the admissible quadruple. 

The following Examples \ref{lip}, \ref{C101n} and \ref{C1T} are exhibited in \cite{hoacta}.
\begin{example}\label{lip}
Let $(K,d)$ be a compact metric space and $Y$ a compact Hausdorff space. Let $0<\alpha\le 1$. Suppose that $B$ is a closed subalgebra of $\Lip((K,d^\alpha))$ which contains the constants and separates the points of $K$, where $d^\alpha$ is the H\"older metric induced by $d$. 
For a metric $d(\cdot,\cdot)$ on $E$, the H\"older metric is defined by $d^\alpha$ for $0<\alpha< 1$. Then 
$\Lip_\alpha((K,d), E)$ is isometrically isomorphic to $\Lip((K,d^\alpha),E)$. 
Suppose that $\widetilde{B}$ is a closed subalgebra of $\Lip((K,d^\alpha),C(Y))$ which contains the constants and separates the points of $K\times Y$. Suppose that $B$ and $\wb$ are self-adjoint. Suppose that 
\[
B\otimes C(Y)\subset \widetilde{B}.
\] 
Let $\mathfrak{M}$ be the Stone-\v Cech compactification of $\{(x,x')\in K^2:x\ne x'\}\times Y$. For $F\in \wb$, let $D(F)$ be the continuous extension to $\mathfrak{M}$ of the function $(F(x,y)-F(x',y))/d^\alpha(x,x')$ on $\{(x,x')\in K^2:x\ne x'\}\times Y$. Then $D:\wb\to C(\mathfrak{M})$ is well defined. 
We have $\|D(F)\|_\infty=L_\alpha(F)$ for every $F\in \wb$.
Hence $\wb$ is a natural $C(Y)$-valuezation of $B$.

There are two typical example of $\wb$ above. 
The algebra $\Lip((K,d^\alpha),C(Y))$ is one. 
The algebras $\Lip((K,d^\alpha))$ and $\Lip((K,d^\alpha),C(Y))$ are self-adjoint (see \cite[Corollary 3]{hoacta}). The inclusions 
\[
\Lip((K,d^\alpha))\otimes C(Y)\subset \Lip((K,d^\alpha),C(Y))
\]
is obvious. Another example of a natural $C(Y)$-valuezation is
$\lip(K,C(Y))$ for $0<\alpha <1$. 
 In fact $\lip(K)$ (resp. $\lip(K,C(Y))$) is a closed subalgebra of $\Lip((K,d^\alpha))$ (resp. $\Lip((K,d^\alpha),C(Y))$ which contains the constants. In this case   Corollary 3 in \cite{hoacta} asserts that $\lip(K)$ separates the points of $K$. As $\lip(K)\otimes C(Y)\subset \lip(K,C(Y))$ we see that $\wb=\lip(K,C(Y))$ separates the points of $K\times Y$. By Corollary 3 in \cite{hoacta}  $\lip(K)$ and $\lip(K,C(Y))$ are self-adjoint. The inclusions 
\[
\lip(K)\otimes C(Y)\subset \lip(K,C(Y))
\]
is obvious.
\end{example}
\begin{example}\label{C101n}
Let $Y$ be a compact Hausdorff space. Then $C^1([0,1],C(Y))$
is a natural $C(Y)$-valuezation of $C^1([0,1])$, where the norm of $f\in C^1([0,1])$ is defined by $\|f\|=\|f\|_\infty+\|f'\|_\infty$ and the norm of $F\in C^{1}([0,1],C(Y))$ is defined by $\|F\|=\|F\|_\infty+\|F'\|_\infty$. It is easy to see that $C^1([0,1])\otimes C(Y)\subset C^1([0,1],C(Y))$.
Let $\mathfrak{M}=[0,1]\times Y$ and $D:C^1([0,1],C(Y))\to C(\mathfrak{M})$ be defined by $D(F)(x,y)=F'(x,y)$ for $F\in C^1([0,1],C(Y))$. Then  $\|F'\|_\infty=\|D(F)\|_{\infty}$ for $F\in C^1([0,1],C(Y))$. 
\end{example}
\begin{example}\label{C1T}
Let $Y$ be a compact Hausdorff space. Then 
$C^1(\mathbb{T},C(Y)))$
is a natural $C(Y)$-valuezation of $C^1({\mathbb {T}})$, where the norm of $f\in C^1(\mathbb{T})$ is defined by $\|f\|=\|f\|_\infty+\|f'\|_\infty$ and the norm of $F\in C^{1}(\mathbb{T},C(Y))$ is defined by $\|F\|=\|F\|_\infty+\|F'\|_\infty$. It is easy to see that $C^1(\mathbb{T})\otimes C(Y)\subset C^1(\mathbb{T},C(Y))$.
Let $\mathfrak{M}=\mathbb{T}\times Y$ and $D:C^1(\mathbb{T},C(Y))\to C(\mathfrak{M})$ be defined by $D(F)(x,y)=F'(x,y)$ for $F\in C^1(\mathbb{T},C(Y))$. Then $\|F'\|_\infty=\|D(F)\|_{\infty}$ for $F\in C^1(\mathbb{T},C(Y))$. 
\end{example}

\section{Isometries on natural $C(Y)$-valuezations}
The following theorem is exhibited in \cite[Theorem 14]{hatori}. It slightly generalize  a similar result for admissible quadruples of type L \cite[Theorem 8]{hoacta}. The proof of Theorem 8 in \cite{hoacta} applies Proposition 3.2 and the following comments in  \cite{hot}. Instead of this we can prove Theorem \ref{main} by Lumer's method, with which a proof is simpler than one given in \cite[Theorem 8]{hoacta}. Refer the detailed proof of Theorm 14 in \cite{hatori} for a proof of Theorem \ref{main}.
\begin{theorem}[Theorem 14 in \cite{hatori}]\label{main}
Suppose that $\wbj$ is a
natural $C(Y_j)$-valuezation of $B\subset C(X_j)$ for $j=1,2$. We assume that 
\[
\|(1\otimes h)F\|_{\wbj}=\|F\|_{\wbj}
\]
for every $F\in \wbj$ and $h\in C(Y_j)$ with $|h|=1$ on $Y_j$ for $j=1,2$. Suppose that $U:\wbo\to \wbt$ is a surjective complex-linear isometry. Then there exists $h\in C(Y_2)$ such that $|h|=1$ on $Y_2$, a continuous map $\varphi:X_2\times Y_2\to X_1$ such that $\varphi(\cdot,y):X_2\to X_1$ is a homeomorphism for each $y\in Y_2$, and a homeomorphism $\tau:Y_2\to Y_1$ which satisfy
\[
U(F)(x,y)=h(y)F(\varphi(x,y),\tau(y)),\qquad (x,y)\in X_2\times Y_2
\]
for every $F\in \wbo$.
\end{theorem}
The weighted composition operator which appears in Theorem \ref{main} has a peculiar form in the sense that the second variable of the composition part depends only on the second variable.  
A composition operator induced by such a homeomorphism is said to be of type BJ in \cite{ho,hot} after the study of 
Botelho and Jamison \cite{bjRocky}. 
\section{The form of $U(1_{\wb_1})$}
Throughout this section we assume that $\wbj$ is a
natural $C(Y_j)$-valuezation of $B\subset C(X_j)$ for $j=1,2$ and 
that $U:\wbo \to \wbt$ is a surjective complex-linear isometry. We assume that $X_2$ is not a singleton in this section. Our main purpose in this section is to show an essence of the proof of  Proposition \ref{absolute value 1}, which is a crucial part of proof of Theorem \ref{main}. Similar proposition and lemmata for admissible quadruples of type L are proved in \cite{hoacta}. Although $\wbj$ in this paper need not be an admissible quadruple of type L, proofs for Proposition \ref{absolute value 1} and Lemmata \ref{1} and \ref{2} are completely the same as that in \cite{hoacta}. Please refer proofs in \cite{hoacta}.
\begin{prop}\label{absolute value 1}
There exists $h\in C(Y_2)$ with $|h|=1$ on $Y_2$ such that $U(1_{\wb_1})=1_{B_2}\otimes h$.
\end{prop}
To prove Proposition \ref{absolute value 1} we apply 
Lemma \ref{2}.
To state Lemma \ref{2} we first define an isometry
 from $\widetilde{B_j}$ into a uniformly closed space of complex-valued continuous functions. 
 Let $j=1, 2$. Define a map
\[
I_j:\wbj \to C(X_j\times Y_j\times \M_j\times \mathbb{T})
\]
by $I_j(F)(x,y,m,\gamma)=F(x,y)+\gamma D_j(F)(m)$ for $F\in \wbj$ and $(x,y,m,\gamma)\in X_j\times Y_j\times M_j\times \mathbb{T}$,where $\mathbb{T}$ is the unit circle in the complex plane. For simplicity we just write $I$ and $D$ instead of $I_j$ and $D_j$ respectively. Scince $D$ is a complex linear map, so is $I$. Put $S_j=X_j\times Y_j\times M_j\times \mathbb{T}$.  For every $F\in \wbj$ the supremum norm $\|I(F)\|_{\infty}$ on $S_j$ of $I(F)$ is written as  
\begin{equation*}
\begin{split}
\|I(F)\|_\infty 
& =\sup\{|F(x,y)+\gamma D(F)(m)|:(x,y,m,\gamma)\in S_j\}\\
& =\sup\{|F(x,y)|:(x,y)\in X_j\times Y_j\}\\
&\qquad 
+\sup\{|D(F)(m)|:m\in \M_j\}\\
&=\|F\|_{\infty(X_j\times Y_j)}+\|D(F)\|_{\infty(\M)}.
\end{split}
\end{equation*}
The second equality holds since $\gamma$ runs through the whole $\mathbb{T}$. 
Therefore we have 
\[
\|I(F)\|_{\infty}=\|F\|_{\infty}
+\|D(F)\|_{\infty}=
 \|F\|_{\wbj}
\]
for every $F\in \wbj$. Since $0= \| D(1)\|_{\infty}$, we have $D(1)=0$ and  $I(1)=1$. Hence $I$ is a complex-linear isometry with $I(1)=1$. In particular, $I(\wbj)$ is a complex-linear closed subspace of $C(S_j)$ which contains $1$. In general $I(\wbj)$ needs not separate the points of $S_j$.

By the definition of the Choquet boundary  $\ch I(\wbt)$ of $I(\wbt)$ (see \cite{ph}), we see that a point $p=(x,y,m,\gamma)\in X_2\times Y_2\times\mathfrak{M}\times \mathbb{T}$ is in $\ch I(\wbt)$ if the point evaluation $\phi_p$ at $p$ is an extreme point of the state space, or equivalently $\phi_p$ is an extreme point of the closed unit ball $(I(\wbt))^*_1$ of the dual space $(I(\wbt))^*$ of $I(\wbt)$.
\begin{lemma}\label{1}
Suppose that 
$(x_0,y_0)\in X_2\times Y_2$ and ${\mathfrak U}$ is an open neighborhood of $(x_0,y_0)$. Then there exists functions $b_0\in B_2$ and $f_0\in C(Y_2)$ such that $0\le F_0\le 1=F_0(x_0,y_0)$ on $X_2\times Y_2$ and $F_0<1/2$ on $X_2\times Y_2\setminus \mathfrak{U}$, where $F_0=b_0\otimes f_0$. Furthermore there exists a point $(x_c,y_c,m_c,\gamma_c)$ in the Choquet boundary for $I_2(\wbt)$ such that $(x_c,y_c)\in {\mathfrak U}\cap F_0^{-1}(1)$ and $\gamma_cD(F_0)(m_c)=\|D(F_0)\|_\infty\ne 0$. 
\end{lemma}
Note that $\gamma_c=1$ if $D(F_0)(m_c)>0$ and $\gamma_c=-1$ if $D(F_0)(m_c)<0$.
\begin{lemma}\label{2}
Suppose that 
$(x_0,y_0)\in X_2\times Y_2$ and ${\mathfrak U}$ is an open neighborhood of $(x_0,y_0)$. Let $F_0=b_0\otimes f_0\in \wbt$ be a function such that $0\le F_0\le 1=F_0(x_0,y_0)$ on $X_2\times Y_2$, and $F_0<1/2$ on $X_2\times Y_2\setminus \mathfrak{U}$. Let $(x_c,y_c,m_c,\gamma_c)$ be a point in the Choquet boundary for $I_2(\wbt)$ such that $(x_c,y_c)\in {\mathfrak U}\cap F_0^{-1}(1)$ and $\gamma_cD(F_0)(m_c)=\|D(F_0)\|_\infty\ne 0$. (Such functions and a point $(x_c,y_c,m_c,\gamma_c)$ exist by Lemma \ref{1}.)
 Then for any $0<\theta<\pi/2$, 
$c_\theta=(x_c,y_c,m_c,e^{i\theta}\gamma_c)$ is also in the Choquet boundary for $I(\wbt)$. 
\end{lemma}
By Lemma \ref{2} we can prove Proposition \ref{absolute value 1} in the same way as the proof of Proposition 9 in \cite{hoacta}.

\section{An application of Lumer's method for a proof of Theorem \ref{main}}
To find isometries Lumer \cite{lum} introduced a useful method which is now called Lumer's method. It involves the notion of Hermitian operators and the fact that $UHU^{-1}$ must be Hermitian if $H$ is Hermitian and $U$ is a surjective isometry. Hermitian operators are usually defined in the notions of the semi inner product. We define it in an equivalent form. A Hermitian element is defined for a unital Banach algebra.
\begin{definition}
Let $\mathfrak{A}$ be a unital Banach algebra. We say that $e\in \mathfrak{A}$ is a Hermitian element if
\[
\|\exp(ite)\|_{\mathfrak{A}}=1
\]
for every $t\in \mathbb{R}$. The set of all Hermitian element of $\mathfrak{A}$ is denoted by $H(\mathfrak{A})$.
\end{definition}
The set of the Hermitian elements $H(M_n(\mathbb{C}))$ in the matrix algebra $M_n({\mathbb C})$ coincides with the set of all Hermitian matrices, and $H(C(Y))=C_{\mathbb R}(Y)$ for the algebra $C(Y)$ of all complex valued continuous functions on a compact Hausdorff space $Y$. In general, for 
a unital $C^*$-algebra $\mathfrak{A}$, the space of all Hermitian elements $H(\mathfrak{A})$ is the space of all self-adjoint elements of $\mathfrak{A}$. A Hermitian element of a unital Banach algebra and a Hermitian operator are sometimes defined in terms of a numerical range, or a semi-inner product. In this paper we define a Hermitian operator by an equivalent form (see \cite{fj1}).
\begin{definition}\label{hermitian}
Let $E$ be a complex Banach space. The Banach algebra of all bounded operators on $E$ is denoted by $B(E)$. We say that $T\in B(E)$ is a Hermitian operator if $T\in H(B(E))$. 
\end{definition}
 The following is a trivial consequence.
\begin{prop}\label{lumer}
Let $E_j$ be a complex Banach space for $j=1,2$. Suppose that $V:E_1\to E_2$ is a surjective isometry and $H:E_1\to E_1$ is a Hermitian operator. Then $VHV^{-1}:E_2\to E_2$ is a Hermitian operator.
\end{prop}
Suppose that $\wb$ is a natural $C(Y)$-valuezation. Since $\nn\cdot\nn:\wb\to \mathbb{R}$ defined by $\nn F\nn=\|D(F)\|_{\infty(\mathfrak{M})}$, $F\in \wb$ is $1$-invariant seminorm (see the definition in \cite{ja}) by the hypothesis on $D:\wb\to C(\mathfrak{M})$. Then $\|F\|_{\wb}=\|F\|_{\infty(X\times Y)}+\|D(F)\|_{\infty(\mathfrak{M})}$, $F\in \wb$ is a natural norm on $\wb$ (see the definition in \cite{ja}). Then Theorem \ref{ja} asserts that a unital surjective isometry $V$ from $\wb$ onto $\wb$ is an isometry from $(\wb, \|\cdot\|_\infty)$ onto $(\wb, \|\cdot\|_\infty)$, too. Hence it is extended to a surjective isometry between the uniform closure of $\wb$ on $X\times Y$. The Stone-Weierstrass theorem asserts that the uniform closure of $\wb$ is $C(X\times Y)$. By the Banach-Stone theorem it is an algebra isomorphism. Hence $V$ is an algebra isomorphism. Hence we have
\begin{prop}\label{iai}
Any surjective unital complex-linear isometry on $\wb$ is an algebra isomorphism.
\end{prop}
Our method of proving Theorem \ref{main} is to find the Hermitian operators. Applying Proposition \ref{iai} we have by Theorem 4 in \cite{ho} that
\begin{prop}\label{hm}
A bounded operator $T$ on $\wb$ is a Hermitian operator if and only if $T(\1)$ is a Hermitian element in $\wb$ and $T=M_{T(\1)}$, the multiplication operator by $T(\1)$. 
\end{prop}
By the similar argument as that in the proof of Proposition 6 in \cite{ho} we have
\begin{prop}\label{hwb}
An element $F\in \wb$ is Hermitian if and only if there exists $u\in C_{\mathbb{R}}(Y)$ such that $F=1\otimes u$.
\end{prop}
Applying propositions above we can prove Theorem \ref{main}. Please refer the proof of Theorem 14 in \cite{hatori}.

\section{Applications of Theorem \ref{main}}
We exhibit applications of Theorem \ref{main}. Corollaries \ref{isoLip},\ref{JPOK},\ref{c101},\ref{c1t} are exhibited in \cite[Section 6]{hoacta}. We omit proofs (see \cite[Section 6]{hoacta}).
\begin{cor}[Corollary 14 in \cite{hoacta}]\label{isoLip}
Let $(X_j,d_j)$ be a compact metric space and $Y_j$ a compact Hausdorff space for $j=1,2$.
Then 
$U:\Lip(X_1,C(Y_1))\to \Lip(X_2,C(Y_2))$ {\rm (}resp. $U:\lip(X_1,C(Y_1))\to \lip(X_2,C(Y_2))${\rm )}
 is a surjective isometry with respect to the norm $\|\cdot\|=\|\cdot\|_\infty+L(\cdot)$ {\rm (}resp. $\|\cdot\|=\|\cdot\|_\infty+L_\alpha(\cdot)${\rm )} if and only if there exists $h\in C(Y_2)$ with $|h|=1$ on $Y_2$, a continuous map $\varphi:X_2\times Y_2\to X_1$ such that $\varphi(\cdot,y):X_2\to X_1$ is a surjective isometry for every $y\in Y_2$, and a homeomorphism $\tau:Y_2\to Y_1$ which satisfy that
\[
U(F)(x,y)=h(y)F(\varphi(x,y),\tau(y)),\qquad (x,y)\in X_2\times Y_2
\]
for every $F\in \Lip(X_1,C(Y_1))$ {\rm (}resp. $F\in \lip(X_1,C(Y_1))${\rm )}.
\end{cor}
Note that if $Y_j$ is a singleton in Corollary \ref{isoLip}, then 
 $\Lip(X_j,C(Y_j))$ {\rm (}resp. $\lip(X_j,C(Y_j))${\rm )} is naturally identified with $\Lip(X_j)$ {\rm (}resp. $\lip(X_j)${\rm )}. 
In this case we have  Example 8 of \cite{jp}.
\begin{cor}[Corollary 15 in \cite{hoacta}]\cite[Example 8]{jp} \label{JPOK}
The map $U:\Lip(X_1)\to \Lip(X_2)$ {\rm (}resp. $U:\lip(X_1)\to \lip(X_2)${\rm )} is a surjective isometry with respect to the norm $\|\cdot\|=\|\cdot\|_{\infty}+L(\cdot)$ {\rm (}resp. $\|\cdot\|=\|\cdot\|_{\infty}+L_\alpha(\cdot)${\rm )} if and only if there exists a complex number $c$ with the unit modulus and a surjective isometry $\varphi:X_2\to X_1$ such that
\[
U(F)(x)=cF(\varphi(x)), \qquad x\in X_2
\]
for every $F\in \Lip(X_1)$ {\rm (}resp. $F\in \lip(X_1)${\rm )}.
\end{cor}
\begin{cor}[Corollary 18 in \cite{hoacta}]\label{c101}
Let $Y_j$ be a compact Hausdorff space for $j=1,2$. The norm $\|F\|$ of $F\in C^{1}([0,1],C(Y_j))$ is defined by $\|F\|=\|F\|_\infty+\|F'\|_\infty$. Then 
$U:C^1([0,1], C(Y_1))\to C^1([0,1],C(Y_2))$ is  a surjective isometry if and only if 
there exists $h\in C(Y_2)$ such that $|h|=1$ on $Y_2$, a continuous map $\varphi:[0,1]\times Y_2\to [0,1]$ such that for each $y\in Y_2$ we have $\varphi(x,y)=x$ for every $x\in [0,1]$ or $\varphi(x,y)=1-x$ for every $x\in [0,1]$, and a homeomorphism $\tau:Y_2\to Y_1$ which satisfy that
\[
U(F)(x,y)=h(y)F(\varphi(x,y),\tau(y)),\qquad (x,y)\in [0,1]\times Y_2
\]
for every $F\in C^1([0,1],C(Y_1))$.
\end{cor}
Note that if $Y_j$ is a singleton in Corollary \ref{c101}, then $C^1([0,1],C(Y_j))$ is $C^1([0,1],{\mathbb C})$. The corresponding result on isometries was given by Rao and Roy \cite{rr}.
\begin{cor}[Corollary 19 in \cite{hoacta}]\label{c1t}
Let $Y_j$ be a compact Hausdorff space for $j=1,2$. The norm $\|F\|$ of $F\in C^{1}(\mathbb{T},C(Y_j))$ is defined by $\|F\|=\|F\|_\infty+\|F'\|_\infty$. Suppose that 
$U:C^1(\mathbb{T}, C(Y_1))\to C^1(\mathbb{T},C(Y_2))$ is a surjective isometry if and only if 
there exists $h\in C(Y_2)$ such that $|h|=1$ on $Y_2$, a continuous map $\varphi:\mathbb{T}\times Y_2\to \mathbb{T}$ and a continuous map $u:Y_2\to \mathbb{T}$ such that for every $y\in Y_2$ $\varphi(z,y)=u(y)z$ for every $z\in \mathbb{T}$ or $\varphi(z,y)=u(y)\bar{z}$ for every $z\in \mathbb{T}$, and a homeomorphism $\tau:Y_2\to Y_1$ which satisfy that
\[
U(F)(z,y)=h(y)F(\varphi(z,y),\tau(y)),\qquad (z,y)\in \mathbb{T}\times Y_2
\]
for every $F\in C^1(\mathbb{T},C(Y_1))$.
\end{cor}
\subsection*{Acknowledgements}
This work was supported by JSPS KAKENHI Grant Numbers JP16K05172, JP15K04921. This work was supported by the Research Institute for Mathematical Sciences, a Joint
Usage/Research Center located in Kyoto University

\end{document}